\documentclass[11pt, reqno]{amsart}
\usepackage{graphics}
\usepackage{amsfonts}

\def\Box{\vcenter{\vbox{\hrule\hbox{\vrule
     \vbox to 8.8pt{\hbox to 10pt{}\vfill}\vrule}\hrule}}}

\newcommand{\Ff}{{\mathbb F}}
\newcommand{\Zz}{{\mathbb Z}}
\newcommand{\Cc}{{\mathbb C}}
\newcommand{\Qq}{{\mathbb Q}}

\newcommand{\gf}{ {{\mathbb F}} }

\def\Tr{\operatorname{Tr}}

\def\PG{\operatorname{PG}}
\def\Sp{\operatorname{Sp}}

\newtheorem{thm}{Theorem}[section]
\newtheorem{lemma}[thm]{Lemma}
\newtheorem{cor}[thm]{Corollary}
\newtheorem{conj}[thm]{Conjecture}
\newtheorem{prop}[thm]{Proposition}

\newtheorem{remark}[thm]{Remark}
\def\Tr{\operatorname{Tr}}

\numberwithin{equation}{section}

\begin{document}

\title[Difference Sets from Symplectic Spreads]
{Skew Hadamard Difference Sets from the Ree-Tits Slice Symplectic
Spreads in $\PG(3,3^{2h+1})$}

\author[Ding, Wang, and Xiang]
{Cunsheng Ding, Zeying Wang, Qing Xiang$^1$}
\thanks{$^1$Research supported in part by NSF Grant DMS 0400411.}

\address{Department of Computer Science,
Hong Kong University of Science and Technology, Clear Water Bay,
Kowloon, Hong Kong, email: {\tt cding@ust.hk}}

\address{Department of Mathematical Sciences, University of Delaware, Newark, DE 19716, USA,
email: {\tt wangz@math.udel.edu}}

\address{Department of Mathematical Sciences, University of Delaware, Newark, DE 19716, USA, email: {\tt
xiang@math.udel.edu}}

\keywords{Difference set, Gauss sum, permutation polynomial,
Ree-Tits slice spread, skew Hadamard difference set, symplectic
spread, Twin prime power difference set}

\date{}

\begin{abstract}
Using a class of permutation polynomials of $\Ff_{3^{2h+1}}$
obtained from the Ree-Tits slice symplectic spreads in
$\PG(3,3^{2h+1})$, we construct a family of skew Hadamard difference
sets in the additive group of $\Ff_{3^{2h+1}}$. With the help of a computer,
we show that these skew Hadamard difference sets are new when $h=2$ and
$h=3$. We conjecture that they are always new when $h>3$.
Furthermore, we present a variation of the classical construction of
the twin prime power difference sets, and show that inequivalent
skew Hadamard difference sets lead to inequivalent difference sets
with twin prime power parameters.
\end{abstract}
\maketitle

\section{Introduction}\label{intro}

Let $G$ be a finite group of order $v$ (written multiplicatively). 
A $k$-element subset $D$ of $G$ is called a {\em $(v,k,\lambda)$ 
difference set\/} if the list of ``differences'' $xy^{-1}$, $x,y\in D$, 
$x\neq y$, represents each nonidentity element in $G$ exactly $\lambda$ times.
As an example of difference sets, we mention the classical Paley
difference set in $(\Ff_q,+)$ consisting of the nonzero squares of 
$\Ff_q$, where $\Ff_q$ is the finite field of order $q$, and $q$ is a prime power congruent to 3 modulo 4. 
Difference sets are the subject of much study in the past 50 years.
We assume that the reader is familiar with the basic theory of
difference sets as can be found in \cite{Baum71}, \cite{lander}, and
\cite[Chap. 6]{bjl}. For a recent survey, see \cite{Xiang}.

A difference set $D$ in a finite group $G$ is called {\it skew
Hadamard} if $G$ is the disjoint union of $D$, $D^{(-1)}$, and
$\{1\}$, where $D^{(-1)}=\{d^{-1}\mid d\in D\}$. The aforementioned
Paley difference set in $(\Ff_q,+)$ is an example of skew Hadamard
difference sets. Let $D$ be a $(v,k,\lambda)$ skew Hadamard
difference set in an abelian group $G$. Then we have
$$1\notin D,\; k=\frac {v-1} {2},\;\mbox {and}\; \lambda=\frac {v-3} {4}.$$
If we employ group ring notation, then in $\Zz[G]$, we have
\begin{eqnarray*}
DD^{(-1)}&=&\frac {v+1} {4}+\frac {v-3} {4}G\\
D+D^{(-1)}&=&G-1,
\end{eqnarray*}
where $D^{(-1)}=\sum_{d\in D}d^{-1}$. Applying any non-principal
(complex) character $\phi$ of $G$ to the above two equations, one has
\begin{equation}\label{charval}
\phi(D)=\frac {-1\pm \sqrt {-v}}{2}.
\end{equation}
Therefore the complex character values of a $(v,k,\lambda)$ skew Hadamard abelian difference set all lie in the quadratic extension $\Qq(\sqrt{-v})$ of $\Qq$. This property of abelian skew difference sets places severe restrictions on these difference sets. Skew
Hadamard difference sets were studied by Johnsen \cite{ecj}, Camion
and Mann \cite{cm}, Jungnickel \cite{jungskew}, and Chen, Xiang and
Seghal \cite{cxs}. The results in \cite{ecj, cm, cxs} can be
summarized as follows:

\begin{thm}\label{skewrestric}
Let $D$ be a $(v,k,\lambda)$ skew Hadamard difference set in an
abelian group $G$. Then $v$ is equal to a prime power $p^m\equiv 3$
(mod 4), and the quadratic residues modulo $v$ are multipliers of
$D$. Moreover, if $G$ has exponent $p^s$ with $s\geq 2$, then $s\leq
(m+1)/4$. In particular, if $v=p^3$ or $p^5$, then $G$ must be
elementary abelian.
\end{thm}

It was conjectured that if an abelian group $G$ contains a skew
Hadamard difference set, then $G$ has to be elementary abelian. This conjecture is still open in general. Theorem~\ref{skewrestric} contains all known results on this conjecture. It was further conjectured some
time ago that the Paley difference sets are the only
examples of skew Hadamard difference sets in abelian groups. This
latter conjecture was recently disproved by Ding and Yuan \cite{DY},
who constructed new skew Hadamard difference sets in $(\Ff_{3^{2h+1}},
+)$ by using certain planar functions related to Dickson
polynomials.

In this paper we construct new skew Hadamard difference sets by
using certain permutation polynomials \cite{bz} from the Ree-Tits
slice symplectic spreads in $\PG(3, 3^{2h+1})$ . While the
construction itself is quite simple (see Section 3), the proof that
the candidate sets are indeed difference sets is not so easy: we had
to resort to a lemma in \cite{cxs} and use Gauss sums and
Stickelberger's theorem on the prime ideal factorization of Gauss
sums. To make the paper self-contained, we include a brief
introduction to Gauss sums here.

Let $p$ be a prime, $q=p^m$. Let $\xi_{p}$ be a fixed complex
primitive $p$th root of unity and let $\Tr_{q/p}$ be the trace from
$\Ff_{q}$ to $\Zz/p\Zz$. Define
$$ \psi: \Ff_{q} \rightarrow \Cc^{*}, \quad \psi(x)=\xi_{p}^{\Tr_{q/p}(x)},$$
which is easily seen to be a nontrivial character of the additive
group of $\Ff_{q}$. Let
$$\chi:\Ff_{q}^* \rightarrow \Cc^{*} $$
be a character of $\Ff_{q}^{*}$ (the cyclic multiplicative group of
$\Ff_q$). We define the {\it Gauss sum} by
$$ g(\chi)=\sum_{a \in \Ff_{q}^*} \chi(a)\psi(a).$$
Note that if $\chi_0$ is the trivial multiplicative character of
$\Ff_q$, then $g(\chi_0)=-1$. Gauss sums can be viewed as the
Fourier coefficients in the Fourier expansion of $\psi|_{\Ff_q^*}$
in terms of the multiplicative characters of
 $\Ff_q$. That is, for every $c\in \Ff_q^*$,
\begin{equation}\label{inv}
\psi(c)=\frac{1} {q-1}\sum_{\chi\in X}g(\chi)\chi^{-1}(c),
\end{equation}
where $X$ denotes the character group of $\Ff_q^*$.

One of the elementary properties of Gauss sums is
\cite[Theorem~1.1.4]{be1}
\begin{equation}\label{eq2.2}
g(\chi)\overline{g(\chi)}=q,\hspace{0.1in}{\rm if}\hspace{0.1in}
\chi\neq \chi_0.
\end{equation}

A deeper result on Gauss sums is Stickelberger's theorem
(Theorem~\ref{stick} below) on the prime ideal factorization of
Gauss sums. We first introduce some notation. Let $a$ be any integer
not divisible by $q-1$. We use $L(a)$ to denote the least positive
integer congruent to $a$ modulo $q-1$. Write $L(a)$ to the base $p$
so that
$$L(a)=a_0+a_1p+\cdots +a_{m-1}p^{m-1},$$
where $0\leq a_i\leq p-1$ for all $i$, $0\leq i\leq m-1$. We define
the {\it digit sum} of $a$ (mod $q-1$) as
$$s(a)=a_0+a_1+\cdots +a_{m-1}.$$
For integers $a$ divisible by $q-1$, we define $s(a)=0$.

Next let $\xi_{q-1}$ be a complex primitive $(q-1)$th root of unity.
Fix any prime ideal $\mathfrak{p}$ in $\Zz[\xi_{q-1}]$ lying over
$p$. Then $\Zz[\xi_{q-1}]/\mathfrak{p}$ is a finite field of order
$q$, which we identify with $\Ff_q$. Let $\omega_{\mathfrak{p}}$ be
the Teichm\"uller character on $\Ff_q$, i.e., an isomorphism
$$\omega_{\mathfrak{p}}: \Ff_q^{*}\rightarrow
\{1,\xi_{q-1},\xi_{q-1}^2,\dots ,\xi_{q-1}^{q-2}\}$$ satisfying
\begin{equation}\label{eq2.3}
\omega_{\mathfrak{p}}(\alpha)\quad ({\rm
mod}\hspace{0.1in}{\mathfrak{p}})=\alpha,
\end{equation}
for all $\alpha$ in $\Ff_q^*$. The Teichm\"uller character
$\omega_{\mathfrak{p}}$ has order $q-1$; hence it generates all
multiplicative characters of $\Ff_q$.

Let $\mathfrak{P}$ be the prime ideal of $\Zz[\xi_{q-1},\xi_p]$
lying above $\mathfrak{p}$. For an integer $a$, let $\nu
_{\mathfrak{P}}(g(\omega_{\mathfrak p}^{-a}))$ denote the
$\mathfrak{P}$-adic valuation of $g(\omega_{\mathfrak p}^{-a})$. The
following classical theorem is due to Stickelberger (see
\cite[p.~7]{lang}, \cite[p.~344]{be1}).

 \begin{thm}\label{stick}
Let $p$ be a prime, and $q=p^m$. Let $a$ be any integer not
divisible by $q-1$. Then
$$\nu_{\mathfrak{P}}(g(\omega_{\mathfrak p}^{-a}))=s(a).$$
\end{thm}

The paper is organized as follows. In Section~\ref{symspread}, we give a brief introduction to symplectic spreads in $\PG(3,q)$, and recall a theorem of
Ball and Zieve \cite{bz} which shows that symplectic spreads in
$\PG(3,q)$ give rise to permutation polynomials of $\Ff_q$ and vice
versa. In particular, we recall a class of permutation polynomials
$f_a(x)$ of $\Ff_{3^m}$, $a\in \Ff_{3^m}$, coming from the Ree-Tits
slice symplectic spreads. In Section~\ref{construction},  we use the
aforementioned permutation polynomials $f_a(x)$ to construct skew
Hadamard difference sets in $(\Ff_{3^m}, +)$. In
Section~\ref{inequiv}, we address the inequivalence issues for skew
Hadamard difference sets in $(\Ff_{3^m},+)$. Finally in
Section~\ref{twinprime}, we present a variation of the classical
construction of the twin prime power difference sets. Also we show
that inequivalent skew Hadamard difference sets can give rise to
inequivalent difference sets with twin prime power parameters.

\section{A Class of Permutation polynomials from the Ree-Tits slice symplectic spreads in $\PG(3,3^{2h+1})$}\label{symspread}

Let $\PG(3,q)$ denote the 3-dimensional projective space over
$\Ff_q$, and let $V=\Ff_q^4$ be the underlying vector space of
$\PG(3,q)$. A {\em spread} of $\PG(3,q)$ is a partition of the
points of the space into lines. Now we equip $V$ with a
non-degenerate alternating form $B: V\times V\rightarrow \Ff_q$. A
spread of $\PG(3,q)$ is called {\em symplectic} if every line of the
spread is totally isotropic with respect to $B$. Since all
non-degenerate alternating forms on $V$ are equivalent, we may assume that
$B$ is defined as follows:
\begin{equation}\label{form}
B((x_0,x_1,x_2,x_3),(y_0,y_1,y_2,y_3))=x_0 y_3-x_3 y_0-x_1 y_2+y_1
x_2
\end{equation}
Then a symplectic spread is a partition of the points of $\PG(3,q)$
into lines such that $B(P,Q)=0$ for any points $P$, $Q$ lying on the
same line of the spread. For readers who are familiar with classical
generalized quadrangles, a symplectic spread of $\PG(3,q)$ is
nothing but a spread of the classical generalized quadrangle
$W_3(q)$. By the Klein correspondence (see \cite{hirschfeld}), a
spread of $W_3(q)$ corresponds to an ovoid of the classical
generalized quadrangle $Q(4,q)$.

In \cite{bz}, it was shown that every symplectic spread of $\PG(3,
q)$ gives rise to a certain family of permutation polynomials of
$\Ff_q$ and vice-visa. Since the symplectic group $\Sp(V)$ leaving
the alternating form in (\ref{form}) invariant acts transitively on
the set of totally isotropic lines, we may assume that the
symplectic spread under consideration contains the line
$$\ell_{\infty}=\langle(0,0,0,1),(0,0,1,0)\rangle.$$

\begin{thm} {\em (\cite{bz})}\label{Ball}
The set of totally isotropic lines
\begin{equation}\label{spreadform}
\ell_{\infty} \cup \{ \langle(0, 1, x, y),(1, 0, -y, g(x,y)\rangle
\mid x, y \in \Ff_q \}
\end{equation}
is a symplectic spread of $\PG(3,q)$ if and only if
$$ x \mapsto g(x,ax-b)+a^2 x $$
is a permutation of $\Ff_q$ for all $a$, $b \in \Ff_q$.
\end{thm}

Table 1 in \cite{bz} lists all known symplectic spreads of
$\PG(3,q)$. For our purpose of constructing new skew Hadamard
difference sets, we are interested in the Ree-Tits slice symplectic
spread, which is a spread having the form
(\ref{spreadform}), with
$$g(x,y)=-x^{2\alpha+3}-y^{\alpha},$$
where $q=3^{2h+1}$ and $\alpha=\sqrt{3q}$. This spread was
discovered by Kantor \cite{kantor} as an ovoid of $Q(4,q)$, which is
a slice of the Ree-Tits ovoid of $Q(6,q)$.

By Theorem~\ref{Ball} the Ree-Tits example gives us a class of
permutation polynomials, namely, the polynomials
$f_a(x)=b^{\alpha}-(g(x,ax-b)+a^2 x)$, $a\in\Ff_q$. Explicitly, we
have
\begin{equation}\label{defpoly}
f_a(x)=x^{2 \alpha+3}+(ax)^{\alpha}-a^2 x.
\end{equation}
As commented in \cite{bz}, the polynomial $f_a$ is remarkable in
that it is a permutation polynomial of $\Ff_q$ whose degree is
approximately $\sqrt{q}$. There are only a handful of known
permutation polynomials with such a low degree. A direct proof that
$f_a(x)$ is a permutation polynomial can be found in \cite{bz}. 

We comment that by going through Table 1 in \cite{bz}, one can see that all other permutation polynomials arising from known symplectic spreads of $\PG(3,q)$, $q$ odd, are linearized permutation polynomials of $\Ff_q$, which will not lead to new skew Hadamard difference sets by the construction described below. That is the reason why we only choose to work with the polynomials $f_a(x)$ defined in (\ref{defpoly}).

\section{A Construction of Skew Hadamard difference
sets}\label{construction}

Throughout this section, $q=3^m$, where $m=2h+1$, $h\geq 0$. For any
$a\in \Ff_q$, let $f_a(x)$ be the polynomial defined in
(\ref{defpoly}). As seen in Section~\ref{symspread}, $f_a(x)$ is a
permutation polynomial of $\Ff_{q}$. For any nonzero $a\in \Ff_q$,
let
\begin{equation}\label{defdiff}
 D_a= \{ f_a(x^2) \mid x\in \Ff_q^*\},
\end{equation}
where $\Ff_q^*=\Ff_q\setminus\{0\}$. We will show that $D_a$ is a skew Hadamard difference set in $(\Ff_q, +)$. We start with the following

\begin{lemma}\label{skew}
For any nonzero $a\in \Ff_q$, we have
$$D_a\cap (-D_a)=\emptyset,$$
and
$$D_a\cup (-D_a) \cup \{ 0\}=\Ff_{q}.$$
\end{lemma}

\begin{proof} Assume that $f_a(x^2)=-f_a(y^2)$ for some $x,y\in
\Ff_q^*$. Then
$$f_a(x^2)=f_a(-y^2).$$
Since $f_a(x)$ is a permutation polynomial of $\Ff_q$, we have
$x^2=-y^2$, which implies that $-1$ is a square in $\Ff_q$. But $-1$
is not a square in $\Ff_q$, since $q=3^m$ and $m$ is odd. Therefore
we reached a contradiction. Hence $D_a\cap (-D_a)=\emptyset$.

Next, clearly we have $f_a(0)=0$. Since $f_a(x)$ is a permutation
polynomial of $\Ff_q$, we see that $f_a(x^2)=0$ if and only if $x=0$
. Therefore $0\not\in D_a$. The second assertion of the lemma now
follows easily. This completes the proof. \end{proof}

We will use the character sum approach (see, e.g.
\cite[p.~318]{bjl}) to prove that $D_a$ is a difference set. Using
this approach, in order to show that $D_a$ is a difference set, we
must prove that for any nontrivial additive character $\psi$ of
$\Ff_q$,
\begin{equation}\label{modulusequ}
\psi(D_a)\overline{\psi(D_a)}=\frac {q+1}{4}.
\end{equation}
It seems difficult to prove directly that (\ref{modulusequ}) holds
for every nontrivial additive characters $\psi$ of $\Ff_q$. We will
use a lemma in \cite{cxs} to bypass this difficulty.

\begin{lemma}{\em (\cite{cxs})}\label{congruence}
Let $G$ be a (multiplicative) abelian $p$-group of order $p^m$,
where $p$ is a prime congruent to $3$ modulo 4, and $m$ is an odd
integer. Let $D$ be a subset of $G$ such that in $\Zz[G]$,
$$D+D^{(-1)}=G -1,$$
and $D^{(t)}=D$ for every nonzero quadratic residue $t$ modulo $p$.
If for every nontrivial character $\phi$ of $G$,
$$\phi(D)\equiv \frac {p^{(m-1)/2}-1} {2}\; \;({\rm mod}\;
p^{(m-1)/2}),$$ then $D$ is a difference set in $G$.
\end{lemma}

The idea of Lemma~\ref{congruence} is that sometimes congruence
properties of $\phi(D)$ can be used to determine the (complex)
absolute value of $\phi(D)$. The proof of the lemma relies on
Fourier inversions, and can be found in \cite{cxs}.

We now state the main theorem of this section.

\begin{thm}\label{main}
Let $a\in \Ff_q^*$, and let $D_a$ be defined as in
(\ref{defdiff}). Then $D_a$ is a skew Hadamard difference set in
$(\Ff_q, +)$.
\end{thm}

\begin{proof}
By Lemma~\ref{skew}, we know that $D_a$ is skew. Since
$1\in \Zz/3\Zz$ is the only nonzero quadratic residue modulo 3, we
certainly have $D_a^{(t)}=D_a$ for every nonzero quadratic residue
$t$ modulo 3. Therefore by Lemma~\ref{congruence}, it suffices to
show that for every nontrivial additive character $\psi_{\beta}:
\Ff_q\rightarrow \Cc^*$,
\begin{equation}\label{additivecong}
\psi_{\beta}(D_a)\equiv \frac {3^{(m-1)/2}-1} {2}\;\; ({\rm mod}\;
3^{(m-1)/2}),
\end{equation}
where $\psi_{\beta}(x)=\xi_3^{\Tr(\beta x)}$, $\xi_3=e^{2\pi i/3}$,
and $\Tr$ is the absolute trace from $\Ff_q$ to $\Ff_3$.

We now compute $\psi_{\beta}(D_a)$. Let $\chi$ be the
(multiplicative) quadratic character of $\Ff_q$. Then

\begin{eqnarray*}
\psi_{\beta}(D_a)&=&\sum_{x\in \Ff_q^*}\psi_{\beta}(f_a(x))\frac {(\chi(x) +1)} {2}\\
&=&\frac {1}{2} (\sum_{x\in \Ff_q^*}\psi_{\beta}(f_a(x))\chi(x) + \sum_{x\in \Ff_q^*}\psi_{\beta}(f_a(x)))\\
&=& \frac {1}{2} (\sum_{x\in \Ff_q^*}\psi_{\beta}(f_a(x))\chi(x)-1),
\end{eqnarray*}
where in the last equality we used the facts that $f_a(x)$ is a
permutation polynomial of $\Ff_q$ and $f_a(0)=0$. From this last
expression for $\psi_{\beta}(D_a)$, we see that (\ref{additivecong})
is equivalent to
\begin{equation}\label{likegausscong}
\sum_{x\in \Ff_q^*}\psi_{\beta}(f_a(x))\chi(x)\equiv 0 \; \; ({\rm
mod} \; 3^{h})
\end{equation}

Let $S_{\beta}=\sum_{x\in \Ff_q^*}\psi_{\beta}(f_a(x))\chi(x)$. We
have
\begin{eqnarray*}
S_{\beta} &=& \sum_{x\in\Ff_q^*}\xi_3^{{\Tr}(\beta x^{2\alpha +3}+(\beta a^{\alpha} -\beta^{\alpha}a^{2\alpha})x^{\alpha})}\chi(x)\\
&=&\sum_{y\in \Ff_q^*}\xi_3^{{\Tr}(\beta y^{\alpha +2} +(\beta a^{\alpha} -\beta^{\alpha}a^{2\alpha})y)}\chi(y)\\
&=&\pm \sum_{y\in \Ff_q^*}\xi_3^{{\Tr}(y^{\alpha +2} +(\beta^{\alpha
-1}a^{\alpha} -\beta^{2\alpha-2}a^{2\alpha})y)}\chi(y)
\end{eqnarray*}

Let $\gamma_a=\beta^{\alpha -1}a^{\alpha}-\beta^{2\alpha
-2}a^{2\alpha}$. If $\gamma_a=0$, then $S_{\beta}$ is a quadratic
Gauss sum, which can be evaluated exactly (see \cite[p.~199]{lidl}).
Indeed, if $\gamma_a=0$, then we have
\begin{eqnarray*}
S_{\beta} &=& \pm \sum_{y \in \Ff_{q}^{*}} \xi_{3}^{\Tr(y^{\alpha+2})}\chi(y) \\
          &=& \pm \sum_{z \in \Ff_{q}^{*}} \xi_{3}^{\Tr(z)}\chi(z) \\
          &=& \pm g(\chi)=\pm \sqrt{-q} =\pm 3^{h}\sqrt{-3} \equiv 0 \pmod {3^h}
\end{eqnarray*}
Hence in this case, $(\ref{likegausscong})$ is true. To finish the
proof, it suffices to prove that when $\gamma_a\neq 0$,

\begin{equation}\label{finalcong}
\sum_{y\in \Ff_q^*}\xi_3^{{\Tr}(y^{\alpha +2}
+\gamma_{a}y)}\chi(y)\equiv0\; \; ({\rm mod} \;3^h)
\end{equation}

Now using Fourier inversion (e.g., see $(\ref{inv})$), we have for
any $y\in \Ff_q^*$,
$$\xi_{3}^{\Tr(y)}=\frac{1}{q-1}\sum_{b=0}^{q-2}g(\omega^{-b})\omega^{b}(y),$$
where $\omega$ is the Teichm\"{u}ller character on $\Ff_{q}$. Then
\begin{eqnarray*}
  &&\sum_{y \in \Ff_{q}^{*}}\xi_{3}^{\Tr(y^{\alpha+2}+\gamma_{a}y)}\chi(y) \\
 &=&\sum_{y \in \Ff_{q}^{*}}\xi_{3}^{\Tr(\gamma_{a}y)}\chi(y)\cdot \frac{1}{q-1}\sum_{b=0}^{q-2}g(\omega^{-b})
\omega^{b}(y^{\alpha+2}) \\
 &=&\sum_{y \in \Ff_{q}^{*}}\xi_{3}^{\Tr(\gamma_{a}y)}\omega^{-\frac{q-1}{2}}(y) \cdot \frac{1}{q-1}
\sum_{b=0}^{q-2}g(\omega^{-b})\omega^{b(\alpha +2)}(y) \\
 &=&\frac{1}{q-1}\sum_{b=0}^{q-2}g(\omega^{-b})\sum_{y \in \Ff_{q}^{*}}\xi_{3}^{\Tr(\gamma_{a}y)}
\omega^{-\frac{q-1}{2}+b(\alpha+2)}(y) \\
 &=&\frac{1}{q-1}\sum_{b=0}^{q-2}g(\omega^{-b})g(\omega^{-\frac{q-1}{2}+b(\alpha+2)})
\omega^{-\frac{q-1}{2}+b(\alpha+2)}(\gamma_{a}^{-1})
\end{eqnarray*}
Hence, we have
\begin{equation}\label{sumofprods}
S_{\beta}=\pm
\frac{1}{q-1}\sum_{b=0}^{q-2}g(\omega^{-b})g(\omega^{-\frac{q-1}{2}+b(\alpha+2)})
\omega^{-\frac{q-1}{2}+b(\alpha+2)}(\gamma_{a}^{-1})
\end{equation}

Fix any prime ideal $\mathfrak{p}$ in $\Zz[\xi_{q-1}]$ lying over
$3$. Let $\mathfrak{P}$ be the prime ideal of $\Zz[\xi_{q-1},\xi_3]$
lying above $\mathfrak{p}$. Since $\nu _{\mathfrak{P}}(3)=2$, we see
that
$$ S_{\beta} \equiv 0 \pmod {3^h} \iff
\nu_{\mathfrak{P}}(S_{\beta}) \ge 2h .$$ Using the expression in
(\ref{sumofprods}) for $S_{\beta}$, we have
\begin{equation}\label{finalequiv}
S_{\beta} \equiv 0 \pmod {3^h} \iff
\nu_{\mathfrak{P}}\left(\sum_{b=0}^{q-2}g(\omega^{-b})
g(\omega^{-\frac{q-1}{2}+b(\alpha+2)})\omega^{-\frac{q-1}{2}+b(\alpha+2)}(\gamma_{a}^{-1})\right)
\ge 2h .
\end{equation}
By Theorem~\ref{stick} and the fact that $g(\chi_0)=-1$, where
$\chi_0$ is the trivial multiplicative character of $\Ff_q$, we have
for any $b$, $0\leq b\leq q-2$,
$$\nu_{\mathfrak{P}}\left(g(\omega^{-b})g(\omega^{-\frac{q-1}{2}+b(\alpha+2)})\right) =
s(b)+s\left (\frac{q-1}{2}-b(\alpha+2)\right ).$$ Therefore if we
can prove that for each $b$, $0 \le b \le q-2$,
\begin{equation}\label{wtinequ}
s(b)+s\left (\frac{q-1}{2}-b(\alpha+2)\right ) \ge 2h,
\end{equation}
then (\ref{finalcong}) will follow. This is exactly what we will do.
In fact, we prove a slightly stronger inequality in
Theorem~\ref{thmappendix}. (Since the proof of
Theorem~\ref{thmappendix} is some what lengthy, we put it in the
Appendix.) Now combine Theorem~\ref{thmappendix} and
Lemma~\ref{congruence}, the proof of the theorem is complete.
\end{proof}

It is of interest to record the following corollary of
Theorem~\ref{main}.

\begin{cor}
Let $q=3^m$, $m=2h+1$, and $\alpha=3^{h+1}$. For any $\beta\in
\Ff_q^*$ and $a\in \Ff_q^*$, we have
$$\sum_{x\in \Ff_q^*}\chi(x)\xi_3^{\Tr(x^{\alpha +2}+(\beta^{\alpha
-1}a^{\alpha}-\beta^{2(\alpha -1)}a^{2\alpha})x)}=\pm\sqrt{-q}.$$
\end{cor}

\section{Inequivalence of skew Hadamard difference
sets}\label{inequiv}

Let $D_1$ and $D_2$ be two $(v,k,\lambda)$ difference sets in an
abelian group $G$. We say that $D_1$ and $D_2$ are {\em
equivalent\/} if there exists an automorphism $\sigma$ of $G$ and an
element $g\in G$ such that $\sigma(D_1)=D_2g$. In this section, we
discuss the inequivalence issues for skew Hadamard difference sets.

\subsection{The known families of skew Hadamard difference sets}

Let $a\in \Ff_q$ and let $n$ be a positive integer. We define the
{\it Dickson polynomial} ${\mathcal D}_n(x,a)$ over $\Ff_q$ by
$${\mathcal D}_n(x,a)=\sum_{j=0}^{\lfloor n/2\rfloor}\frac {n}{n-j}{n-j\choose j}(-a)^jx^{n-2j},$$
where $\lfloor n/2\rfloor$ is the largest integer $\leq n/2$. It is
well known that the Dickson polynomial ${\mathcal D}_n(x,a)$, $a\in
\Ff_q^*$, is a permutation polynomial of $\Ff_q$ if and only if
$\gcd(n,q^2-1)=1$ (see \cite[p.~356]{lidl}). Let $m$ be a positive
odd integer. For any $u\in \Ff_{3^m}^*$, define
$$g_u(x)={\mathcal D}_5(x^2, -u)=x^{10}-ux^6-u^2x^2.$$
It was proved in \cite{DY} that when $m$ is a positive odd integer
and $u\in \Ff_{3^m}^*$, ${\rm Image}(g_u)\setminus \{0\}$ is a skew
Hadamard difference set in $(\Ff_{3^m}, +)$. For convenience, we set
$$DY(u)=\{x^{10}-ux^6-u^2x^2 \mid x\in \Ff_{3^m}^*\},$$ 
and call these {\it the Ding-Yuan difference sets}. We have the following

\begin{prop}\label{prop-knownDS}
All previously known skew Hadamard difference sets are equivalent to one of the
following:
\begin{enumerate}
\item The Paley difference set $P$ in $\gf_q$, where $q \equiv 3 \pmod{4}$
      is a prime power.
\item The Ding-Yuan difference set $DY(1)$ in $\gf_{3^m}$,
      where $m$ is odd.
\item The Ding-Yuan difference set $DY(-1)$ in $\gf_{3^m}$,
      where $m$ is odd.
\end{enumerate}
\end{prop}

\begin{proof}
First of all, it can be checked directly that ${\mathcal D}_5(-x, u)
= -{\mathcal D}_5(x, u)$ and
\begin{eqnarray}\label{eqn-Dick2}
b^5 {\mathcal D}_5(x, a) = {\mathcal D}_5(bx, b^2a), \; \forall a, b
\in \Ff_q
\end{eqnarray}
Setting $a=-1$ in (\ref{eqn-Dick2}), we have
$$
b^5 {\mathcal D}_5(x^2, -1) = {\mathcal D}_5(bx^2, -b^2).
$$
Thus, we have $DY(b^2)=b^5 DY(1)$ if $b$ is a nonzero square in $\Ff_{3^m}$; and $
DY(b^2)=- b^5 DY(1)$ if $b$ is a nonsquare. Hence for any nonzero square $u\in
\Ff_{3^m}$, $DY(u)$ is equivalent to $DY(1)$.

Similarly, we can prove that for any nonsquare $u\in \Ff_{3^m}$, $DY(u)$ is equivalent to $DY(-1)$.

Combining the above observation with the fact that the Paley family
and the Ding-Yuan family were the only previously known skew
Hadamard difference sets, we see that the proof of the proposition
is complete.
\end{proof}

With the help of a computer, it was verified in \cite{DY} that the
three skew Hadamard difference sets $P$, $DY(1)$ and $DY(-1)$ in $(\gf_{3^m}, +)$ are
all equivalent when $m=3$, but they are indeed pairwise inequivalent
when $m=5$ and $7$. It is very likely that the three difference sets
$P$, $DY(1)$ and $DY(-1)$ are pairwise inequivalent for all odd $m>7$, although this is
not proved rigorously.

\subsection{The inequivalence issues for the difference sets $D_a$}
 We now turn to the difference sets $D_a$ constructed in Section 3. First we prove the following

\begin{prop}
Let $m=2h+1$ be a positive integer and let $a\in \Ff_{3^m}^*$. The
skew Hadamard difference sets $D_a$ in $(\Ff_{3^m}, +)$ constructed
in Section 3 are equivalent to one of the following:\

(1) The difference set $D_{1}$ in $(\Ff_{3^m}, +)$.\

(2) The difference set $D_{-1}$ in $(\Ff_{3^m}, +)$.\\
\end{prop}

\begin{proof} Using the definition of $f_a(x)$ in (\ref{defpoly}), it can
be checked that
$$b^{2\alpha+3}f_{a}(\frac{x}{b})=f_{ab^{\alpha+1}}(x), \; \forall  b\in \Ff_{3^m}^{*}.$$

Assume that $a$ is a nonzero square in $\Ff_{3^m}$. Since
$\gcd(\alpha+1,q-1)=2$, one can find $ \zeta \in \Ff_{3^m}^{*}$ such
that
$$a \zeta^{\alpha+1}=1.$$
Hence
\begin{equation} \label{trans}
   \zeta^{2 \alpha+3}f_{a}(\frac{x^2}{\zeta})=f_{1}(x^2).
\end{equation}
We note that if  $\zeta$ is a square, then $
\{f_{a}(\frac{x^2}{\zeta})\mid x\in \Ff_{3^m}^{*} \}=
\{f_{a}(x^2)\mid x\in \Ff_{3^m}^{*}\}=D_{a}$; and if $\zeta$ is a
nonsquare, then $\{f_{a}(\frac{x^2}{\zeta})\mid x\in \Ff_{3^m}^{*}
\}= \{f_{a}(-x^2)\mid x\in \Ff_{3^m}^{*}\}=\{-f_{a}(x^2)\mid x\in
\Ff_{q}^{*}\}=-D_{a}$. Therefore
\begin{equation}\label{newzeta}
\{ \zeta^{2 \alpha+3} f_a(\frac{x^2}{\zeta})\mid x \in \Ff_{3^m}^{*}
\}=\zeta^{2 \alpha+3}D_a, \; {\rm or}\; -\zeta^{2 \alpha+3}D_a .
\end{equation}
Combining (\ref{newzeta}) with (\ref{trans}), we see that
$D_{a}$ is equivalent to $D_{1}$.

Similarly, we can show that $D_{a}$ is equivalent to $D_{-1}$ when
$a$ is a nonsquare in $\Ff_{3^m}$.
\end{proof}

Since equivalent difference sets give rise to isomorphic symmetric
designs, which have the same $p$-rank and Smith normal form, we may
use $p$-ranks and Smith normal forms to distinguish inequivalent
difference sets. See \cite{Xiang} for a recent survey of results on
this subject. Unfortunately, skew Hadamard difference sets with the
same parameters have the same $p$-rank \cite[pp. 297--299]{Jung1}
and the same Smith normal form \cite{MW98}. Thus in order to
distinguish inequivalent skew Hadamard difference sets, we have to
use some other techniques.

It seems not easy to settle completely the question whether the
difference sets $D_1$ and $D_{-1}$ are inequivalent to the
previously known families stated in Proposition~\ref{prop-knownDS}.
With the aid of a computer, we will show that the skew Hardamard
difference sets $D_1$ and $D_{-1}$ in $(\Ff_{3^m},+)$ are new when
$m=5$ and $7$. (We mention that when $m=3$, the difference sets $D_1$ and $D_{-1}$ are equivalent to the Paley difference set in $\Ff_{3^3}$.)

Let $D$ be a difference set in $(\gf_q, +)$. For any 2-subset $\{a, b\}\subset \Ff_q^*$, we define
$$
T\{a,b\}: = |D \cap (D+a) \cap (D+b)|.
$$
These numbers $T\{a,b\}$ are called the {\em triple intersection numbers}, which were used to distinguish inequivalent difference sets in 1971 by Baumert \cite[p.144]{Baum71}.

We shall use the triple intersection numbers to distinguish the skew difference sets of this paper 
from the earlier ones in the cases where $m=5$ and $m=7$. We use $P$ and $RT(a)$
to denote the Paley difference set and the difference set $D_a$ from
Section~\ref{construction}, respectively.

With the help of Magma~\cite{cannon}, the maximum and minimum triple
intersection numbers of these difference sets in $\gf_{3^7}$ are
computed and listed below:
\begin{eqnarray*}
\begin{array}{lcccc}
\mbox{Difference set} & & \mbox{ Minimum (when $m=7$)} &  & \mbox{ Maximum (when $m=7$)} \\
P                    & &  261                   &  & 284 \\
DY(1)                 & &  246                   &  & 300    \\
DY(-1)                & &  248                   &  & 297    \\
RT(1)                  & &  250                   &  & 295    \\
RT(-1)                 & &  249                   &  & 296
\end{array}
\end{eqnarray*}
Hence the five difference sets are pairwise inequivalent when $m=7$.
It then follows from Proposition \ref{prop-knownDS} that the skew
difference sets $RT(1)$ and $RT(-1)$ are new when $m=7$.

When $m=5$, the maximum and minimum triple intersection numbers of these difference
sets in $\gf_{3^m}$ are computed and listed below:
\begin{eqnarray*}
\begin{array}{lcccc}
\mbox{Difference set} & & \mbox{ Minimum (when $m=5$)} &  & \mbox{ Maximum (when $m=5$)} \\
P                   & &  26                   &  & 33 \\
DY(1)                 & &  23                  &  & 36    \\
DY(-1)                & &  24                   &  & 35    \\
RT(1)                & &  24                   &  & 35    \\
RT(-1)                & &  24                   &  & 35
\end{array}
\end{eqnarray*}
In fact, in this case $DY(-1)$, $RT(1)$ and $RT(-1)$ have the same
set of triple intersection numbers, i.e., $\{i: 24 \leq i \leq
35\}$. We further compute the multiplicities of these triple intersection
numbers for these three cases. We find the following data:
\begin{eqnarray*}
\begin{array}{lcccc}
\mbox{Difference set} & & \mbox{Triple intersection numbers with multiplicities ($m=5$)}   \\
DY(-1)                & &  24^{75}25^{435}26^{1155}27^{2385}\cdots 35^{120}\\
RT(1)                & &  24^{75}25^{330}26^{1155}27^{2535}\cdots 35^{105}\\
RT(-1)                & &  24^{90}25^{330}26^{1095}27^{2655}\cdots 35^{120},
\end{array}
\end{eqnarray*}
where the exponents denote multiplicities. Since the multiplicities of the (triple) intersection number 27 are pairwise distinct for the three cases, we conclude that 
$DY(-1)$, $RT(1)$ and $RT(-1)$ are pairwise inequivalent when $m=5$. Hence, the five difference sets $P$, $DY(1)$, $DY(-1)$, $RT(1)$, and $RT(-1)$ are pairwise 
inequivalent when $m=5$. It then follows from Proposition \ref{prop-knownDS} that the
skew difference sets $RT(1)$ and $RT(-1)$ are new when $m=5$.

Based on the above evidence, we make the following conjecture.

\begin{conj}
The five difference sets $P$, $DY(1)$,
$DY(-1)$, $RT(1)$ and $RT(-1)$ in $(\Ff_{3^m},+)$ are pairwise
inequivalent for all odd $m>7$.
\end{conj}

\section{Difference sets with twin prime power
parameters}\label{twinprime}

In this section we present a variation of the classical construction
of the twin prime power difference sets. Using this variation we
will show that inequivalent skew Hadamard difference sets can give
rise to inequivalent difference sets with twin prime power
parameters. We first recall the construction of the twin prime power
difference sets. As usual, we denote the (multiplicative) quadratic
character of a finite field by $\chi$.

\begin{thm}{\em (Stanton and Sprott~\cite{ss})}\label{classical}
Let $q$ and $q+2$ be odd prime powers. Then the set
$$D=\{(x,y)\mid x\in \Ff_q^*, \;y\in\Ff_{q+2}^*, \;
\chi(x)=\chi(y)\}\cup \{(x,0)\mid x\in\Ff_q\}$$ is a
$(4n-1,2n-1,n-1)$ difference set in $(\Ff_q,+)\times (\Ff_{q+2},+)$,
where $n=\frac {(q+1)^2}{4}$.
\end{thm}

For a proof of Theorem~\ref{classical}, we refer the reader to
\cite{ss} or \cite[p.~354]{bjl}. For convenience, we will refer the
parameters $(4n-1,2n-1,n-1)$, $n=\frac {(q+1)^2}{4}$, $q$ an odd
prime power, as the twin prime power parameters. We now give a
variation of the above construction.

\begin{thm}\label{vary}
Let $q$ and $q+2$ be prime powers, and let $q\equiv 3$ (mod 4). Let
$E$ be a skew Hadamard difference set in $(\Ff_q,+)$. Then the set
$$D=\{(x,y)\mid x\in E, \;y\in\Ff_{q+2}^*, \; \chi(y)=1\}\cup
\{(x,y)\mid x\in -E, \;y\in\Ff_{q+2}^*, \; \chi(y)=-1\}\cup
\{(x,0)\mid x\in\Ff_q\}$$ is a $(4n-1,2n-1,n-1)$ difference set in
$(\Ff_q,+)\times (\Ff_{q+2},+)$, where $n=\frac {(q+1)^2}{4}$.
\end{thm}

Noting that the nontrivial character values of a skew Hadamard
difference set are given by (\ref{charval}), one can easily give a
character theoretic proof for Theorem~\ref{vary}. We leave this to the reader as an exercise.

\begin{remark}
(1). We remark that if $q$ and $q+2$ are both prime powers, and
$q\equiv 1$ (mod 4), then we can similarly use a skew Hadamard
difference set in $\Ff_{q+2}$ to construct a difference set in
$(\Ff_q,+)\times (\Ff_{q+2},+)$ with twin prime power parameters.

(2). One further generalization of Theorem~\ref{vary} goes as
follows. With the assumptions in Theorem~\ref{vary}, let $Q$ be any
$(q+2,\frac {q+1}{2}, \frac {q-3}{4}, \frac {q+1}{4})$ partial
difference set in $(\Ff_{q+2}, +)$, $0\not\in Q$. (See
\cite[p.~230]{bjl} for the defintion of partial difference set.)
Then the set
$$D'=\{(x,y)\mid x\in E, \;y\in Q\}\cup
\{(x,y)\mid x\in -E, \;y\in\Ff_{q+2}^*\setminus Q\}\cup \{(x,0)\mid
x\in\Ff_q\}$$ is a $(4n-1,2n-1,n-1)$ difference set in
$(\Ff_q,+)\times (\Ff_{q+2},+)$, where $n=\frac {(q+1)^2}{4}$.

\end{remark}

In view of the fact that there exist inequivalent skew Hadamard
difference sets in $(\Ff_q,+)$, the following theorem is of
interest.

\begin{thm}\label{inequivtwin}
Let $q$ and $q+2$ be prime powers, and let $q\equiv 3$ (mod 4). Let
$E$ and $F$ be inequivalent skew Hadamard difference sets in
$(\Ff_q,+)$. Then the two difference sets
$$D=\{(x,y)\mid x\in E, \;y\in\Ff_{q+2}^*, \; \chi(y)=1\}\cup
\{(x,y)\mid x\in -E, \;y\in\Ff_{q+2}^*, \; \chi(y)=-1\}\cup
\{(x,0)\mid x\in\Ff_q\}$$ and
$$D'=\{(x,y)\mid x\in F, \;y\in\Ff_{q+2}^*, \; \chi(y)=1\}\cup
\{(x,y)\mid x\in -F, \;y\in\Ff_{q+2}^*, \; \chi(y)=-1\}\cup
\{(x,0)\mid x\in\Ff_q\}$$ are inequivalent.
\end{thm}

\begin{proof}
Assume that $D$ and $D'$ are equivalent difference sets in
$G=(\Ff_q,+)\times (\Ff_{q+2},+)$. Then there exists an automorphism
$\alpha$ of $G$ and an element $(b_1,b_2)\in G$ such that
\begin{equation}\label{concrete}
\alpha(D)=D' + (b_1,b_2).\end{equation}
We will show that $E$ and $F$ are equivalent.

For convenience, we define
$$A_1=\{(x,y)\mid x\in E, \;y\in\Ff_{q+2}^*, \; \chi(y)=1\}\cup
\{(x,y)\mid x\in -E, \;y\in\Ff_{q+2}^*, \; \chi(y)=-1\},$$
$$A_2=\{(x,y)\mid x\in F, \;y\in\Ff_{q+2}^*, \; \chi(y)=1\}\cup
\{(x,y)\mid x\in -F, \;y\in\Ff_{q+2}^*, \; \chi(y)=-1\},$$ and
$$B=\{(x,0)\mid x\in\Ff_q\}.$$
So $D=A_1\cup B$, $D'=A_2\cup B$, and (\ref{concrete}) can be
written as
\begin{equation}\label{moreconcrete}
\alpha(A_1)\cup \alpha(B)=(A_2+ (b_1,b_2))\cup (B +
(b_1,b_2)).\end{equation}

Since $\gcd(q,q+2)=1$, we have ${\rm Aut}(G)\cong {\rm
Aut}(\Ff_q,+)\times {\rm Aut}(\Ff_{q+2},+)$. Hence there exist $f\in
{\rm Aut}(\Ff_q,+)$ and $g\in {\rm Aut}(\Ff_{q+2},+)$ such that
$\alpha(x,y)=(f(x),g(y))$ for all $(x,y)\in G$.

We claim that $b_2=0$. If not, then there exists a $y\in
\Ff_{q+2}^*$ such that $g(y)=b_2$. Note that
$B+(b_1,b_2)=\{(x,b_2)\mid x\in \Ff_q\}$. By (\ref{moreconcrete}),
we must have
$$\{(f(x),g(y))\mid x\in E\}=\{(x,b_2)\mid x\in \Ff_q\},$$
or
$$\{(f(x),g(y))\mid x\in -E\}=\{(x,b_2)\mid x\in \Ff_q\},$$
according as $\chi(y)=1$ or $\chi(y)=-1$. However both equalities
are clearly impossible by comparing the cardinalities of the sets
involved. This proves that $b_2=0$. It follows that
$\alpha(B)=B+(b_1,0)$ and
\begin{equation}\label{finequiv}
\alpha(A_1)=A_2 + (b_1,0).
\end{equation}
Let $y\in \Ff_{q+2}^*$ such that $g(y)=1$. From (\ref{finequiv}), we
see that
$$\{(f(x),g(y))\mid x\in E\}=\{(x+b_1,1)\mid x\in F\},$$
or
$$\{(f(x),g(y))\mid x\in -E\}=\{(x+b_1,1)\mid x\in F\},$$
according to $\chi(y)=1$ or $\chi(y)=-1$. So $f(E)=F+b_1$ or
$-f(E)=F+b_1$. This proves that $E$ and $F$ are equivalent
difference sets in $(\Ff_q,+)$.
\end{proof}

Combining Theorem~\ref{inequivtwin} with the results in
Section~\ref{inequiv}, we see that whenever $3^{2h+1}\pm 2$ ($h>1$) is a
prime power, there exist difference sets with twin prime power
parameters that are inequivalent to the classical twin prime power
difference sets. To indicate that there are $h>1$ such that $3^{2h+1}\pm 2$
are prime powers, we mention the following specific examples: $3^5-2=241$
is a prime, $3^9-2=19681$ is a prime, and $3^{15}+2=14348909$ is also a prime.

\section{Appendix}\label{inequality}

In this appendix, we give the promised proof of (\ref{wtinequ}).
Throughout this section, $m=2h+1$ is a positive odd integer,
$q=3^m$, $r=\frac{m+1}{2}=h+1$, and $\alpha=3^{\frac{m+1}{2}}=3^r$.
Our goal is to prove
\begin{thm}\label{thmappendix}
For each $a$, $0 \le a \le q-2$, we have
\begin{equation}\label{strong}
s(a)+s\left (\frac{q-1}{2}-a(\alpha+2)\right ) \ge m, \end{equation}
where $s(a)$ is the digit sum of $a$ defined in Section 1.
\end{thm}
First of all, we observe that the only $a$, $0 \le a \le q-2$,
satisfying
$$\frac {q-1}{2}-a(\alpha+2)\equiv 0 \; ({\rm mod}\; q-1)$$
is $a=\frac {q-1}{2}$. For $a=\frac {q-1}{2}$, we have $s(a)=m$ and
$s\left (\frac{q-1}{2}-a(\alpha+2)\right )=0$. So certainly
(\ref{strong}) holds for $a=\frac {q-1}{2}$. Therefore in our
discussion below, we will always assume that $a\neq\frac {q-1}{2}$
(and $\frac{q-1}{2}-a(\alpha+2)\not\equiv 0$ (mod $q-1$)).

A sequence $\{u_i\}_{i \in \Zz}$ is called periodic with period $m$
if $u_i=u_j$ whenever $i \equiv j \pmod{m}$. All sequences in this
section are periodic with period $m$. Let $a$ be an integer
satisfying $0\leq a\leq q-2$ and $a\neq \frac {q-1}{2}$. Write
$$a=\sum_{i=0}^{m-1} a_{i}3^{i}, \quad a_i\in\{0,1,2\}$$
and extend $a_0$, $a_1$, $\cdots$, $a_{m-1}$ to a periodic sequence with period $m$.
We have
\begin{eqnarray*}
\frac{q-1}{2}-(3^r+2)a&=&\frac{q-1}{2}-3^{r}a-3a+a\\
                &\equiv &\sum_{i=0}^{m-1}(1-a_{i-r}-a_{i-1}+a_i)3^{i} \pmod {3^m-1} \\
                &\equiv &\sum_{i=0}^{m-1}\left(1+(2-a_{i-r})+(2-a_{i-1})+a_{i} \right)3^{i} \pmod {3^{m}-1}\\
                &=& \sum_{i=0}^{m-1}(5+a_i-a_{i-1}-a_{i-r})3^{i}  \quad
\end{eqnarray*}
For each $i$, let
$$b_i=5+a_i-a_{i-1}-a_{i-r}.$$
It is easily seen that $b_i\in \{1,2,3,\ldots,7\}$. Write
$$\sum_{i=0}^{m-1}b_{i}3^{i}\equiv\sum_{i=0}^{m-1}s_{i}3^{i} \pmod{3^{m}-1}$$
with $s_{i}\in\{0,1,2\}$. By Theorem 13 of \cite{hx} (adapted to the
ternary case), there exists a sequence $\{c_i\}$ such that
\begin{equation}
\forall i , \; s_{i}=b_{i}-3c_{i}+c_{i-1},
\end{equation}
where $c_i\in\{0,1,2,3\}$ is the carry from the $i$th digit to the
$(i+1)$th digit in the modular summation of $\frac{q-1}{2}$,
$-3^ra$, $-3a$ and $a$. Note that
\begin{eqnarray*}
& &s(a)+s\left (\frac{q-1}{2}-(3^{r}+2)a\right )\\
&=&\sum_{i=0}^{m-1}a_{i}+\sum_{i=0}^{m-1}\left((5+a_i-a_{i-1}-a_{i-r})-3c_{i}+c_{i-1}\right)\\
&=&5m-2\sum_{i=0}^{m-1}c_{i}
\end{eqnarray*}
So in order to prove Theorem~\ref{thmappendix}, it suffices to prove
\begin{equation}
\sum_{i=0}^{m-1}c_{i}\le 2m
\end{equation}

Since $\gcd(r,m)=\gcd(r,2r-1)=1$, for any fixed $i$, the sequence
$c_i$, $c_{i-r}$, $c_{i-2r}$, $\ldots$, $c_{i-(m-1)r}$ is a
rearrangement of $c_0$, $c_1$, $\ldots$, $c_{m-1}$. In the
following, we will also frequently use the facts that $2r \equiv 1
\pmod{m}$, $c_{i-1}=c_{i-2r}$, $c_{i-2}=c_{i-4r}$, and so on.


\begin{lemma}\label{base}
If $c_{i}=3$, then $c_{i-1}=2$ and $c_{i-r}\le 2$, $a_{i}=2$,
$a_{i-1}=a_{i-r}=0$.
\end{lemma}

\begin{proof} Note that $s_{i}=b_{i}-3c_{i}+c_{i-1} \ge 0$, $1 \le b_{i} \le
7$, $0 \le c_{i-1} \le 3$. If $c_{i}=3$, then
\begin{equation}\label{3}
6 \le b_{i} \le 7, \quad 2 \le c_{i-1} \le 3.
\end{equation}
Assume to the contrary that $c_{i-1}=3$. Since
$$s_{i-1}=b_{i-1}-3c_{i-1}+c_{i-2} \ge 0 ,$$
we have
\begin{equation}\label{4}
 6 \le b_{i-1} \le 7, \quad 2 \le c_{i-2} \le 3 .
\end{equation}
From the lower bounds on $b_i$ and $b_{i-1}$ in (\ref{3}) and
(\ref{4}), we have
\begin{eqnarray*}
 6 &\leq & b_{i}=5+a_{i}-a_{i-1}-a_{i-r},\; {\rm and}\\
 6 &\leq & b_{i-1}=5+a_{i-1}-a_{i-2}-a_{i-r-1}.
\end{eqnarray*}
Adding up the two inequalities, we get
$$10+a_{i}-a_{i-r}-a_{i-2}-a_{i-r-1} \ge 12,$$
which implies that
$$a_{i}=2,a_{i-r}=a_{i-r-1}=a_{i-2}=0.$$
We use the following table to summarize the above information.
\begin{equation*}
{\rm A}:=\left[ \begin{array}{ccccc} a_{i} & a_{i-r} & a_{i-1} & a_{i-r-1}  & a_{i-2} \\
                        2    & 0       & \ge 0   & 0          & 0
     \end{array} \right]
\end{equation*}
Since
$$ b_i=5+a_{i}-a_{i-1}-a_{i-r} \ge 6 ,$$
using the information in Table A, we have
$$ a_{i-1} \le 1.$$
Since
$$ b_{i-1}=5+a_{i-1}-a_{i-2}-a_{i-r-1} \ge 6 ,$$
again using the information in Table A, we have
$$ a_{i-1} \ge 1.$$
Hence $a_{i-1}=1$. Therefore we can update the entries in Table A as
follows.
\begin{equation*}
{\rm A}=\left[ \begin{array}{ccccc} a_{i} & a_{i-r} & a_{i-1} & a_{i-r-1} & a_{i-2} \\
                        2    &  0      &  1      &  0        &  0
         \end{array} \right]
\end{equation*}
It follows that $b_{i-1}=5+a_{i-1}-a_{i-2}-a_{i-r-1}=6$. Since
$s_{i-1}=b_{i-1}-3c_{i-1}+c_{i-2} \ge 0$ and $c_{i-1}=3$, we have
$c_{i-2}=3$. Combining this with $s_{i-2}=b_{i-2}-3c_{i-2}+c_{i-3}
\ge 0$, we obtain $b_{i-2} \ge 6$.\

Since
\begin{eqnarray*}
     b_{i-1}&=&5+a_{i-1}-a_{i-2}-a_{i-r-1} \ge 6, \\
     b_{i-2}&=&5+a_{i-2}-a_{i-3}-a_{i-r-2} \ge 6,
\end{eqnarray*}
adding up these two inequalities, we get
$$10+a_{i-1}-a_{i-r-1}-a_{i-3}-a_{i-r-2} \ge 12,$$
which implies that
$$a_{i-1}=2,a_{i-r-1}=a_{i-3}=a_{i-r-2}=0. $$
But this is in contradiction with the previous conclusion that
$a_{i-1}=1$ as shown in Table A. Hence $c_{i-1} \ne 3$. By (\ref{3})
we must have $c_{i-1}=2$.

Combining the fact $c_{i-1}=2$, $c_{i}=3$ with
$s_{i}=b_{i}-3c_{i}+c_{i-1} \ge 0$, we have $b_i=7$. Recall that
$$b_{i}=5+a_{i}-a_{i-1}-a_{i-r}\leq 7.$$
We obtain
$$a_{i}=2, a_{i-1}=a_{i-r}=0.$$
Now $b_{i-r}=5+a_{i-r}-a_{i-r-1}-a_{i-1}=5-a_{i-r-1} \le 5$,
$s_{i-r}=b_{i-r}-3c_{i-r}+c_{i-r-1} \ge 0$, and $c_{i-r-1} \le 3$,
we conclude that $c_{i-r} \le 2$. This completes the proof.
\end{proof}

\begin{lemma} \label{base2}
If $c_{i}=3$, $c_{i-r}=c_{i-1}=2$, then $c_{i-r-1} \le 2$. That is,
\begin{equation*}
\left [ \begin{array}{ccc} c_{i} & c_{i-r} & c_{i-1} \\
                             3   &  =2     &  =2
        \end{array} \right]
\implies
\left[ \begin{array}{cccc}c_{i} & c_{i-r} & c_{i-1} & c_{i-r-1}\\
                             3   &  =2     &  =2    &  \le 2
       \end{array} \right].
\end{equation*}
\end{lemma}

\begin{proof} Assume to the contrary that $c_{i-r-1}=3$. By Lemma \ref{base},
we have $c_{i-r-2}=2$, $c_{i-2} \le 2$, $a_{i-r-1}=2$, and
$a_{i-r-2}=a_{i-2}=0$. Since
$$ s_{i-1}=b_{i-1}-3c_{i-1}+c_{i-2} \ge 0, $$
and $c_{i-1}=2$, $c_{i-2} \le 2$, we have $b_{i-1}\ge 4$. By
assumption $c_{i}=3$. It follows from Lemma \ref{base} that
$a_{i}=2$, $a_{i-1}=a_{i-r}=0$. We use the following table to
summarize the above information.
\begin{equation*}
{\rm A}:=\left[ \begin{array}{cccccc} a_{i} & a_{i-r} & a_{i-1} & a_{i-r-1} & a_{i-2} & a_{i-r-2} \\
                                2    &  0      &   0     &   2       &  0      &   0
         \end{array} \right]
\end{equation*}
Recall that
$$b_{i-1}=5+a_{i-1}-a_{i-2}-a_{i-r-1}.$$
Using the information in Table A, we have $b_{i-1}=5+0-0-2=3$,
which contradicts the previous conclusion that $b_{i-1}\ge 4$. This
completes the proof.
\end{proof}

\begin{thm}\label{thmappendix1}
Let $t \ge 3$ be an integer. If $c_{i}=3$,
$c_{i-r}=c_{i-2r}=\cdots=c_{i-tr}=2$, then $a_{i}=2$, $a_{i-r} \le
1$, $a_{i-2r} \le 1$, $\ldots$, $a_{i-(t-1)r} \le 1$,
$a_{i-(t-2)r}+a_{i-(t-1)r} \le 1$, and $c_{i-tr-r} \le 2$.
Furthermore, if $c_i=3$, $c_{i-r}=c_{i-2r}=\cdots=c_{i-tr}=2$ and
also $c_{i-tr-r}=2$, then $a_{i-tr} \le 1$ and
$a_{i-(t-1)r}+a_{i-tr} \le 1$.
\end{thm}

\begin{proof}
We will use induction on $t$. When $t=3$, the assumptions are
$c_{i}=3$, and $c_{i-r}=c_{i-2r}=c_{i-3r}=2$ (i.e.,
$c_{i-r}=c_{i-1}=c_{i-r-1}=2$). We will show that $a_i=2$, $a_{i-r}
\le 1$, $a_{i-2r}=a_{i-1} \le 1$, $a_{i-r}+a_{i-1} \le 1$, and
$c_{i-3r-r}=c_{i-2} \le 2$.

Since $c_i=3$, by Lemma \ref{base}, we have
\begin{equation}\label{5}
a_{i}=2, a_{i-r}=0, a_{i-1}=0.
\end{equation}
It remains to show that $c_{i-2}\le 2$. Assume to the contrary that
$c_{i-2}=3$, by Lemma \ref{base}, we have $c_{i-3}=2$ and $c_{i-2-r}
\le 2$, $a_{i-2}=2$, $a_{i-3}=a_{i-2-r}=0$. We summarize the
information in the following table
\begin{equation*}
{\rm A}:=\left[ \begin{array}{ccccccc} a_i & a_{i-r} & a_{i-1} & a_{i-r-1} & a_{i-2} & a_{i-r-2} & a_{i-3} \\
                                 2  &  0      &  0      &  \le 2    & 2       &    0      &  0
         \end{array} \right].
\end{equation*}
Since
\begin{align*}
 s_{i-r-1}&=b_{i-r-1}-3c_{i-r-1}+c_{i-r-2} \ge 0,\\
 s_{i-1}  &=b_{i-1}-3c_{i-1}+c_{i-2} \ge 0,\\
 c_{i-r-1}&=2,\; c_{i-r-2} \le 2,\; c_{i-1}=2,\; c_{i-2}=3,
\end{align*}
we see that $b_{i-r-1}\ge 4$ and $b_{i-1}\ge 3$. Using the
information in Table A, we find that
$$b_{i-r-1}=5+a_{i-r-1}-a_{i-r-2}-a_{i-2}=3+a_{i-r-1}.$$
So $b_{i-r-1}\ge 4$ implies that $a_{i-r-1} \ge 1$. Again using the
information in Table A, we find that
$$b_{i-1}=5+a_{i-1}-a_{i-2}-a_{i-r-1}=3-a_{i-r-1}.$$
So $b_{i-1}\ge 3$ implies that $a_{i-r-1} \le 0$, which contradicts
with the previous conclusion that $a_{i-r-1} \ge 1$. Therefore we
must have $c_{i-2} \le 2$.

Next we show that if $c_i=3$ and
$c_{i-r}=c_{i-2r}=c_{i-3r}=c_{i-4r}=2$ (i.e.,
$c_{i-r}=c_{i-1}=c_{i-r-1}=c_{i-2}=2$), then $a_{i-3r}=a_{i-r-1} \le
1$, and $a_{i-1}+a_{i-r-1} \le 1$. Note that (\ref{5}) is still
true. Since
\begin{align*}
s_{i-r}&=b_{i-r}-3c_{i-r}+c_{i-r-1} \ge 0 ,\\
s_{i-1}&=b_{i-1}-3c_{i-1}+c_{i-2} \ge 0 ,\\
c_{i-r}&=c_{i-1}=c_{i-r-1}=c_{i-2}=2,
\end{align*}
we have
\begin{eqnarray}
4 &\leq & b_{i-r}=5+a_{i-r}-a_{i-r-1}-a_{i-1}, \label{8}\\
4 &\leq & b_{i-1}=5+a_{i-1}-a_{i-2}-a_{i-r-1}.\label{9}
\end{eqnarray}
Adding up (\ref{8}) and (\ref{9}), we get
$$10+a_{i-r}-2a_{i-r-1}-a_{i-2} \ge 8 .$$
As $a_{i-r}=0$ (see (\ref{5})), the above inequality becomes
\begin{equation*}
2a_{i-r-1}+a_{i-2} \le 2
\end{equation*}
Since $a_{i-2} \ge 0$, we have
$$a_{i-r-1} \le 1 .$$
Noting that $a_{i-1}=0$ (see (\ref{5})), we have
$$a_{i-1}+a_{i-r-1} \le 1.$$
This finishes the proof in the case where $t=3$.

Assume that the theorem is proved for $t=k-1\ge 3$. We will prove
the theorem for $t=k$. So assume that $c_{i}=3$,
$c_{i-r}=c_{i-2r}=\cdots=c_{i-kr}=2$. By induction hypothesis, we
have
\begin{equation}\label{10}
\begin{split}
a_{i}=2, a_{i-r} \le 1, a_{i-2r} \le 1, \ldots, a_{i-(k-2)r} \le 1, \\
a_{i-(k-3)r}+a_{i-(k-2)r} \le 1.
\end{split}
\end{equation}
Since it is also assumed that $c_{i-kr}=2$, we have
\begin{equation}\label{11}
a_{i-(k-1)r} \le 1, \quad a_{i-(k-2)r}+a_{i-(k-1)r} \le 1 .
\end{equation}
Now we show that $c_{i-kr-r} \le 2$. Assume to the contrary that
$c_{i-kr-r}=3$. Then by Lemma \ref{base}, we have $c_{i-kr-r-1}=2$,
$c_{i-kr-1} \le 2$, and $a_{i-kr-r}=2$, $a_{i-kr-r-1}=a_{i-kr-1}=0$.
As before we summarize the information in the following table.
\begin{equation*}
{\rm B}:=\left[ \begin{array}{cccccc} a_{i-(k-2)r} & a_{i-(k-1)r} & a_{i-kr} & a_{i-kr-r} & a_{i-kr-1} & a_{i-kr-r-1} \\
                          \le 1       &   \le 1      &   \ge 0  &    2       &   0        &   0
           \end{array}   \right].
\end{equation*}
Since
\begin{eqnarray*}
s_{i-(k-1)r}&=&b_{i-(k-1)r}-3c_{i-(k-1)r}+c_{i-kr-r} \ge 0, \\
s_{i-kr}    &=&b_{i-kr}-3c_{i-kr}+c_{i-kr-1} \ge 0, \\
c_{i-(k-1)r}&=&c_{i-kr}=2, \; c_{i-kr-r}=3,\; c_{i-kr-1} \le 2;
\end{eqnarray*}
we have
\begin{eqnarray}
3 &\leq & b_{i-(k-1)r}=5+a_{i-(k-1)r}-a_{i-kr-r}-a_{i-kr},\label{22}\\
4 &\leq & b_{i-kr}   =5+a_{i-kr}-a_{i-kr-1}-a_{i-kr-r}.\label{23}
\end{eqnarray}
Adding up (\ref{22}) and (\ref{23}), we get
$$ 10+a_{i-(k-1)r}-2a_{i-kr-r}-a_{i-kr-1} \ge 7 .$$
Since $a_{i-kr-r}=2$, we have
$$ a_{i-(k-1)r}-a_{i-kr-1} \ge 1 ,$$
which implies that $a_{i-(k-1)r} \ge 1$. Combining this with the
information on $a_{i-(k-1)r}$ in Table B, we have
\begin{equation} \label{24}
a_{i-(k-1)r}=1.
\end{equation}
Thus we can update the information in Table B as follows.
\begin{equation*}
{\rm B}=\left[ \begin{array}{cccccc} a_{i-(k-2)r} & a_{i-(k-1)r} & a_{i-kr} & a_{i-kr-r} & a_{i-kr-1} & a_{i-kr-r-1} \\
                          \le 1       &   =1     &   \ge 0  &    2       &   0        &   0
           \end{array}   \right].
\end{equation*}
Using the updated Table B and (\ref{22}) (respectively, (\ref{23})), we get $a_{i-kr} \le 1$ (respectively, $a_{i-kr} \ge 1$).
Hence
\begin{equation} \label{25}
a_{i-kr}=1.
\end{equation}
Since
\begin{eqnarray*}
s_{i-(k-3)r}&=&b_{i-(k-3)r}-3c_{i-(k-3)r}+c_{i-(k-1)r} \ge 0, \\
s_{i-(k-2)r}&=&b_{i-(k-2)r}-3c_{i-(k-2)r}+c_{i-kr} \ge 0 ,\\
c_{i-(k-3)r}&=&c_{i-(k-2)r}=c_{i-(k-1)r}=c_{i-kr}=2 ;
\end{eqnarray*}
we have
\begin{eqnarray}
4 &\leq & b_{i-(k-3)r}=5+a_{i-(k-3)r}-a_{i-(k-1)r}-a_{i-(k-2)r}, \label{26}\\
4 & \leq &
b_{i-(k-2)r}=5+a_{i-(k-2)r}-a_{i-kr}-a_{i-(k-1)r}.\label{27}
\end{eqnarray}
Adding up (\ref{26}) and (\ref{27}), we get
\begin{equation}\label{28}
 10+a_{i-(k-3)r}-2a_{i-(k-1)r}-a_{i-kr} \ge 8 .
\end{equation}
Noting that $a_{i-(k-1)r}=a_{i-kr}=1$, we obtain from (\ref{28}) and
(\ref{27}) that
\begin{equation} \label{29}
a_{i-(k-3)r} \ge 1 , \quad a_{i-(k-2)r} \ge 1 ,
\end{equation}
which implies that
$$a_{i-(k-3)r}+a_{i-(k-2)r} \ge 2 ,$$
contradicting with $ a_{i-(k-3)r}+a_{i-(k-2)r} \le 1$ in (\ref{10}).
Therefore
\begin{equation*}
c_{i-kr-r} \le 2 .
\end{equation*}

Finally, assume that $c_i=3$,
$c_{i-r}=c_{i-2r}=\cdots=c_{i-kr}=c_{i-kr-r}=2$. From the
conditions, we know that (\ref{10}), (\ref{11}), (\ref{26}) and
(\ref{27}) still hold. Since
\begin{eqnarray*}
 s_{i-(k-1)r}&=&b_{i-(k-1)r}-3c_{i-(k-1)r}+c_{i-kr-r} \ge 0 ,\\
 c_{i-(k-1)r}&=&c_{i-kr-r}=2,
\end{eqnarray*}
we have
\begin{equation} \label{30}
4\leq b_{i-(k-1)r}=5+a_{i-(k-1)r}-a_{i-(k+1)r}-a_{i-kr}.
\end{equation}
Adding up (\ref{26}) and (\ref{27}), (\ref{27}) and (\ref{30}),
respectively, we get
\begin{eqnarray}
2a_{i-(k-1)r}+a_{i-kr}-a_{i-(k-3)r} &\le& 2 , \label{31}\\
2a_{i-kr}+a_{i-(k+1)r}-a_{i-(k-2)r} &\le& 2 . \label{32}
\end{eqnarray}
Adding up (\ref{31}) and (\ref{32}), we get
$$ 2a_{i-(k-1)r}+3a_{i-kr} \le 4+(a_{i-(k-3)r}+a_{i-(k-2)r})-a_{i-(k+1)r} .$$
Since $a_{i-(k-3)r}+a_{i-(k-2)r} \le 1$, it follows that
\begin{equation} \label{33}
2a_{i-(k-1)r}+3a_{i-kr} \le 5-a_{i-(k+1)r} \le 5 .
\end{equation}
Now we would like to show that $a_{i-(k-1)r}+a_{i-kr} \le 1$. Assume
to the contrary that $a_{i-(k-1)r}+a_{i-kr} >1$. Since $a_{i-(k-1)r}
\le 1$, by (\ref{33}), we have
\begin{equation} \label{34}
 a_{i-(k-1)r}=1, \quad a_{i-kr}=1.
\end{equation}
Combining (\ref{27}) with (\ref{34}), we get
$$ a_{i-(k-2)r} \ge 1.$$
So $a_{i-(k-2)r}+a_{i-(k-1)r} \ge 2 $, contradicting with
$a_{i-(k-2)r}+a_{i-(k-1)r} \le 1$ in (\ref{11}). Hence
\begin{equation*}
a_{i-(k-1)r}+a_{i-kr} \le 1 ,
\end{equation*}
which in turn implies
$$a_{i-kr} \le 1.$$
This completes the proof.
\end{proof}

\begin{cor}\label{cor1}
If $c_i=3$, then $c_{i-r} \le 2$.  Let $t \ge 1$ be an integer. If
$c_i=3$, $c_{i-r}=c_{i-2r}=\cdots=c_{i-tr}=2$, then $c_{i-tr-r} \le
2$.
\end{cor}

\begin{proof} The first assertion follows directly from Lemma \ref{base}.
For the second assertion, when $t=1$ (respectively, $t=2$), the corollary
follows from Lemma~\ref{base} (respectively, Lemma \ref{base2}).
When $t\geq 3$, the corollary follows directly from Theorem~\ref{thmappendix1}.
\end{proof}

\begin{cor}\label{cor2}
Let $t$ be an integer satisfying $1 \le t \le m$. If
$c_i=c_{i-tr}=3$ and $c_{i-r} \le 2$, $c_{i-2r} \le 2$, $\ldots$,
$c_{i-(t-1)r} \le 2$, then there exists $\ell$, $1 \le \ell \le
{t-1}$ such that $c_{i-\ell r} <2$.
\end{cor}

\begin{proof} Using Corollary \ref{cor1}, we see that the condition $c_{i}=c_{i-tr}=3$ implies that $t\ne 1$. Assume to the contrary that for every $\ell \in \{1,2, \ldots, t-1\}$,
we have $c_{i-\ell r}=2$. That is,
$$c_{i}=3, c_{i-r}=c_{i-2r}=\cdots=c_{i-(t-1)r}=2.$$
Then by Corollary \ref{cor1}, we have $c_{i-tr}\le 2$, contradicting
with our assumption that $c_{i-tr}=3$. This completes the proof.
\end{proof}

\noindent {\bf Proof of Theorem \ref{thmappendix}.} As stated
before, it suffices to prove that
\begin{equation}
\sum_{i=0}^{m-1}c_i \le 2m \label {*}
\end{equation}
If $c_i \le 2$, for all $i$, $ 0 \le i \le {m-1}$, then the
inequality (\ref{*}) of course holds. So we assume that there exists
an $h$, $0\leq h\leq m-1$, such that $c_h=3$. Since $\gcd(m,r)=1$,
we see that the sequence $c_h,c_{h-r},\ldots ,c_{h-(m-1)r}$ is just
a permutation of $c_0, c_1,\ldots , c_{m-1}$. We assume that
$c_{h-i_1 r}=c_{h-i_2 r}=\ldots =c_{h-i_s r}=3$, where $0=i_1 <i_2
<\cdots i_s <m$, $s\geq 1$, and $c_{h-j r}\leq 2$ for each $j\in
\{0,1,\ldots ,m-1\}\setminus \{i_1,i_2,\ldots ,i_s\}$. By
Corollary~\ref{cor1}, we have $i_2\geq i_1+2$, $i_3\geq i_2+2$,
\ldots ,$i_s\geq i_{s-1}+2$, and $m\geq i_s+2$. Using
Corollary~\ref{cor2}, we can bound the sum of the entries in each
segment as follows:
\begin{align*}
c_{h-i_1 r}+c_{h-(i_1 +1)r}+\ldots +c_{h-(i_2 -1)r} &\leq 2(i_2 - i_1),\\
c_{h-i_2 r}+c_{h-(i_2 +1)r}+\ldots +c_{h-(i_3 -1)r} &\leq 2(i_3 - i_2),\\
\vdots \\
c_{h-i_s r}+c_{h-(i_s +1)r}+\ldots +c_{h-(m -1)r} &\leq 2(m - i_s).
\end{align*}
Summing up the above inequalities, we obtain (\ref{*}). The proof of the theorem is complete.

\vspace{0.1in}
\noindent{\bf Acknowledgement.} The authors thank an anonymous referee for his/her helpful comments.

\bibliographystyle{amsalpha}

\end{document}